\documentstyle[11pt,amscd,amstex, epsfig]{article}
\def\al{\alpha}

\def\dim{{\rm dim\; }}

\def\text#1{{\em #1}}

\def\ep{\varepsilon}

\def\be{\begin{equation}}
\def\ee{\end{equation}}
\def\bear{\begin{eqnarray}}
\def\eear{\end{eqnarray}}
\def\best{\begin{eqnarray*}}
\def\eest{\end{eqnarray*}}
\def\pf{{\bf Proof}: }

\renewcommand{\theequation}{\arabic{section}.\arabic{equation}}
\newtheorem{theorem}{Theorem}[section]
\newtheorem{prop}[theorem]{Proposition}
\newtheorem{proposition}[theorem]{Proposition}
\newtheorem{lemma}[theorem]{Lemma}
\newtheorem{cor}[theorem]{Corollary}

\newtheorem{defn}[theorem]{Definition}

\newtheorem{ex}[theorem]{Example}

\def\non{\noindent}

\def\pf{\non {\bf Proof. }}
\def\QED{\nopagebreak \hskip .1in { $\Box$ }\penalty10000 %
\hskip\parfillskip \par  }

\def\Si{\Sigma}

\def\dbar{\overline \partial_J}

\def\Z{{ \Bbb Z}}

\def\P{{ \Bbb P}}

\def\cx{{ \Bbb C}}

\def\w{\omega}

\def\M{{\cal M}}

\def\H{{\cal H}}

\def\PH{{\cal H}}
\def\dbara{\overline\partial_{J_\al}}

\newcommand{\SF}{\overline{{\cal F}}_{g,k}(X,A)}
\newcommand{\F}{\overline{{\cal F}}_{g,k}^{l}(X,A)}

\newcommand{\CM}{\overline{{\cal M}}}

\begin{document}

\date{\empty}
\title{\bf Family  Gromov-Witten Invariants for  K\"{a}hler Surfaces\vskip.2in}
\author{Junho Lee\\
University of Minnesota\\
Minneapolis, MN, 55454}

\addtocounter{section}{0}


\maketitle

\vskip.15in

\begin{abstract}
We use analytic methods to define Family  Gromov-Witten Invariants for
K\"{a}hler surfaces.  We prove that
these are well-defined invariants of the deformation class of the
K\"{a}hler structure.
\end{abstract}

\vskip.4in

Gromov-Witten invariants are counts of holomorphic curves in a
symplectic manifold $X$.  To define them
using the analytic approach one  chooses an  almost  complex
structure $J$ compatible with the symplectic
structure  and considers the set of maps $f:\Sigma\to X$ from  Riemann surfaces
$\Sigma$ which satisfy the (nonlinear elliptic)
$J$-holomorphic map equation
\begin{equation}\label{PSHE}
\overline{\partial}_Jf=0.
\label{1}
\end{equation}
     After compactifying the moduli space of such
maps, one  imposes constraints, requiring, for example, that the
image of the map passes through  specified points.  With the
right number of constraints and a generic $J$, the number of such
maps is finite.  That number is a GW invariant; it depends
only on the symplectic structure of $X$.

\bigskip

There are some beautiful conjectures
about what the counts of holomorphic curves on  K\"{a}hler
surfaces {\em ought to be} (\cite{v},\cite{kp},\cite{yz},\cite{g}).
However, as currently defined,
the GW invariants with those many point constraints corresponding to
the dimension of holomorphic curves on K\"{a}hler surfaces with
$p_{g}\geq 1$ are all zero!
The reason is that the dimension of the GW moduli space
for {\em generic} almost  complex structures
is strictly less than the dimension of the space of relevant
holomorphic curves on K\"{a}hler surfaces.
This discrepancy in dimension occurs because
K\"{a}hler structures  are very special,
namely the linearizations obtained from (\ref{PSHE}) are not onto
for K\"{a}hler structures when $p_{g}\geq 1$
---
Donaldson details this in \cite{d}(see also \cite{bl3}).

Clearly a new
version of the invariants is needed --- one which counts the
relevant holomorphic curves. Work in that direction is just beginning.
Bryan and Leung (\cite{bl1},\cite{bl2}) defined such   invariants for
K3 and Abelian surfaces by using
the Twistor family; they  were also able to calculate their
invariants in important cases.
Behrend-Fantechi \cite{bf} have
defined  invariants for a more general
class of algebraic surfaces using algebraic geometry.
We approach the same issues using the geometric analysis approach to GW
invariants.

\bigskip

Given a K\"{a}hler manifold $(X,\omega,J,g)$ we construct a
$2p_g$-dimensional family of  elements $K_{J}(f,\alpha)$ in
$\Omega^{0,1} (f^*TX)$, where $\alpha$ is a real part of a
holomorphic 2 form. We then  modify the $J$-holomorphic map
equation (1) by considering the pairs $(f,\alpha)$ satisfying
\begin{equation}
\overline{\partial}_Jf= K_{J}(f,\alpha).
\label{eq0.2}
\end{equation}
The solutions of this equation form a  moduli space   whose dimension
is $2p_g$ larger than the dimension of the usual GW
moduli space.

Because $\alpha$  ranges over a vector space
compactness is an issue.  Here things get
interesting because there are instances when the moduli space for
(\ref{eq0.2}) is {\em not} compact.  In
fact, when the map represents a component of a canonical divisor the
moduli space is {\em never} compact.    Nevertheless, there is a
simple analytic criterion --- the
uniform boundedness of the energy of the map and the $L^2$ norm of
$\alpha$ --- that ensures that the
moduli space  is  compact.

\bigskip

\begin{theorem}\label{thm01}
Let $(X, J)$  be a  K\"{a}hler surface and fix a genus $g$ and a
class $A\in H_2(X,\Z)$.
Denote by $C(J)$ the supremum of $E(f)+||\alpha||_{L^2}$
over all $(J,\alpha)$-holomorphic maps
from genus $g$ curves into $X$ which
represent $A$. If $C(J)$ is finite, then
the family
     GW invariants
$$
GW^{J,{\cal H}}_{g,k}(X,A)
$$
are well-defined.
They are invariant under deformations $\{J_{t}\}$ of the K\"{a}hler structure
with $C(J_{t})$ bounded.
Furthermore, if $A$ is a $(1,1)$ class then all the maps
which contribute to these invariants are in fact
$J$-holomorphic.
\end{theorem}

The last sentence of Theorem \ref{thm01}  means that the invariants
for $(1,1)$ classes  are counts of
holomorphic curves in
$(X, J)$.  That is not the same as saying the invariants are
enumerative, since there is no claim that
each curve is counted with multiplicity one.  But it does mean that
the family GW invariants, which {\em
a priori} are counts of maps which are holomorphic with respect to
families of almost complex structures
on $X$, are in fact calculable from the complex geometry of $(X, J)$ alone.

\medskip

Theorem \ref{thm01} yields well-defined family GW invariants provided
there is a finite
energy bound $C(J)$.  Following the Kodaira classification of
surfaces, we verify the energy bound
   case-by-case using geometric arguments.  That yields the following
cases where the family GW
invariants are well-defined.

\begin{prop}\label{0.15}
The moduli space for a class $A$ is compact,
and hence the family GW
invariants $GW^{J,\H}_{g,k}(X,A)$ are well-defined,   when
$(X,J)$ is

\medskip

  (a) \  a K3 or Abelian surface with $A\ne 0$,

\smallskip

(b) \  a minimal  elliptic surface $\pi:E\to C$ with Kodaira
dimension 1 with $-A\cdot (\mbox{fiber
class}) \neq
\mbox{deg }\pi_*(A)$ , and

\smallskip

(c) \   a minimal surface of general type and $A$ is of type $(1,1)$
such that $A$ is not
a linear combination of components of the canonical class $K$.

\end{prop}

Next consider a minimal surface $(X,J)$  of general type.  As
mentioned above,  Proposition~\ref{0.15}
does not apply to classes $A$ which
are linear combinations of  components of the canonical class (\,such
as  $A=mK$ for $m\geq 1$\,) because, for those $A$, the family moduli
space is not compact (\,see Example~\ref{Ex:CPT}\,).
However, we can still define invariants for such classes  by adapting
the approach used by Ionel and Parker ([IP1,IP2]) to define GW
invariants.  That gives  the following extension of Proposition~\ref{3}.

\begin{prop}\label{3}
Let $(X,J)$ be a minimal surface of general type.
Suppose $A$ is of type $(1,1)$ and the class $A$ and
the genus $g$ satisfy
\begin{equation*}
  -K\cdot A\ +\ g-1\ \geq\ 0.
\end{equation*}
Then there are well-defined invariants
for the class $A$ and the genus $g$ which coincide
with the family invariants whenever
$A$ also satisfies Condition (c) of Proposition~\ref{0.15}.
\end{prop}

\vskip 0.8 cm

Section 1 gives the definition of a $(J,\alpha)$-holomorphic map and
some of the analytic consequences
of that definition, most notably an expression for the  energy   in
terms of  pullback of
the symplectic form and the form $\alpha$.
Section 2 begins by
describing  the relation between
a complete linear system $|C|$ --- or  more generally a Severi variety ---
and the moduli space of $(J,\alpha)$-holomorphic maps.  That leads us
to consider the family of
$(J,\alpha)$-holomorphic maps in which $\alpha$ is
the real part of holomorphic 2-form; the corresponding family moduli
space should be an analytic
version of the Severi variety.  As partial justification of that
view, we prove the last statement of
Theorem  \ref{thm01}:   any $(J,\alpha)$-holomorphic map which
represents a (1,1) class is in fact
holomorphic (theorem \ref{L:GT}).

Section 3 summarizes the analytic results which lead to the definition  of the
  family GW-invariants.  That involves constructing the  virtual moduli cycle
by adapting the method  of Li and Tian \cite{LT}.  Thus defined, the
family invariants satisfy
a Composition Law analogous to those of ordinary
GW-invariants.

\smallskip

Section 4 contains examples of K\"{a}hler surfaces
with $p_{g}\geq 1$ with well-defined family invariants.   We focus on
minimal surfaces and establish the
results summarized in Propositions \ref{0.15} and
~\ref{3} above.   For the case of
K3 and Abelian surfaces we prove
that our
  family GW-invariants agree  with the invariants defined by Bryan and
Leung.  That is done in the course
of the  proof of Theorem~\ref{C:K3} by relating the holomorphic
2-forms to the  Twistor family.

\smallskip

The appendix contains a brief discussion of how the   family GW
invariants defined here relate to those
defined by  Behrend and Fantachi in \cite{bf}.

\bigskip
\noindent {\bf Acknowledgments : } I would like to thank most
sincerely my advisor Prof. Thomas Parker for his guidance and
helpful discussions. Without his help this paper would not have
been possible. I am also very grateful to Prof. Ronald Fintushel
for his interest in this work and wish to thank the referee for
comments which helped to improve the presentation of this paper.



\vskip 1cm

\setcounter{equation}{0}
\section{$(J,\alpha)$-holomorphic maps}
\label{section1}
\bigskip

A $J$-holomorphic map into an almost complex manifold $(X,J)$ is a
map $f:\Sigma\to X$ from a complex curve
$\Sigma$  (a closed Riemann surface with complex structure $j$) whose
differential is complex linear.
Equivalently, $f$ is a solution of the $J$-holomorphic map equation
$$
\dbar f=0 \qquad \mbox{ where}\qquad \dbar f = \dfrac{1}{2}(df+Jdfj).
$$
In this section we will
show that when $X$ is
four-dimensional there is natural infinite-dimensional family of
almost complex structures parameterized
the $J$-anti-invariant 2-forms on $X$.

\bigskip

Let $(X,J)$ be a 4-dimensional almost K\"{a}hler manifold
with the hermitian triple $(\omega,J,g)$.
Using $J$, we can decompose  $\Omega^{2}(X)\otimes {\Bbb C}$ as
$
  \Omega^{2}(X)\otimes\,{\Bbb C} =\
  \Omega^{1,1}_{J}\oplus\,
  (\,\Omega^{2,0}_{J}\oplus\Omega^{0,2}_{J}).
$
This leads to the following decomposition
$$
  \Omega(X)\ =\ \Omega^{+}_{J}(X)\ \oplus\ \Omega^{-}_{J}(X)
$$
where $\Omega^{+}_{J}(X)=(\,\Omega^{1,1}_{J})_{\,\Bbb R}$\ and
$\Omega^{-}_{J}(X)=
  (\,\Omega^{2,0}_{J})_{\,\Bbb R}\, (\, =
  (\,\Omega^{2,0}_{J}\oplus\Omega^{0,2}_{J})_{\,\Bbb R}\,).$
Note that $\alpha\in\Omega^{-}_{J}(X)$ iff\ \
$\alpha(Ju,Jv)=-\alpha(u,v)$.

\medskip

\begin{defn}
A 2-form $\alpha$ is called $J$-anti-invariant
if $\alpha$ is in $\Omega^{-}_{J}(X)$.
Each $\alpha$ in $\Omega^{-}_{J}(X)$
defines an endomorphism $K_{\alpha}$ of $TX$ by the equation
\begin{equation}
\langle u,K_{\alpha}v \rangle = \alpha(u,v).
\end{equation}
\label{def1.1}
\end{defn}
It follows that
\bear
\langle K_{\alpha}u , v \rangle =
               - \langle u , K_{\alpha} v \rangle,\ \
J K_{\alpha}=-K_{\alpha} J,\ \ {\rm and}\ \
\langle Ju , K_{\alpha}u \rangle = 0.
\label{1.3}
\eear

\bigskip

\begin{defn}
For $\alpha\in\Omega^{-}_{J}(X)$,
a map $f:\Sigma\to X$ is called ($J,\alpha$)-holomorphic if
\bear
\label{1.defJalphaholo}
\dbar f = K_{J}(f,\alpha)
\eear
where\  $K_{J}(f,\alpha)=K_\al(\partial f\circ j)=
\dfrac{1}{2}K_{\alpha}(df-Jdfj)j$.
\label{def1.2}
\end{defn}

\bigskip

The next proposition and its corollary  list some pointwise relations
involving the quantities that
appear in the ($J,\alpha$)-holomorphic equation.  We state these
first for general $C^1$ maps, then
specialize to  ($J,\alpha$)-holomorphic maps.


\begin{proposition}
\label{P:Basic}
Fix a metric within the conformal class  $j$ and let $dv$ be the
associated volume form.  Then for any
$C^1$ map $f$ we have the pointwise equalities

\bigskip

\begin{tabular}{ll}
(a)  \ \  $|\dbar f|^{2}\ dv = \dfrac{1}{2}\,|df|^{2} \ dv
                 - f^{*}\omega$, \hspace{0.5cm} &
(b) \ \ $\langle \dbar f\,  ,\,  K_{J}(f,\alpha)\rangle \ dv
                  = f^{*}\alpha$, \\   \medskip
(c)\ \   $K_{\alpha}^{2}=-|\alpha|^{2}\,Id$,
&(d)\ \  $|K_{J}(f,\alpha)|^{2}\ dv = f^{*}(\,|\alpha|^{2})
           \left(\dfrac{1}{2}|df|^{2}\ dv + f^{*}\omega\right). $
\end{tabular}
\end{proposition}

\pf Fix a point $p\in \Sigma$ and  an orthonormal
basis $\{ e_{1}, e_{2} = je_{1} \}$ of $T_{p}\Sigma$.  Setting
       $v_{1}=df(e_{1})$ and $v_{2}=df(e_{2})$, we have
$2\dbar f(e_{1}) = v_{1} +Jv_{2}$ and
$2K_{J}(f,\alpha)(e_{1})  = K_{\alpha}v_{2} - JK_{\alpha}v_{1}$, and similarly
       $2\dbar f(e_{2}) = v_{2}-Jv_{1}$ and
$2K_{J}(f,\alpha)(e_{2})   = -K_{\alpha}v_{1} - JK_{\alpha}v_{2}$.
Therefore,
\best
4|\dbar f|^{2}   = |v_{1} +Jv_{2}|^{2} + |v_{2}-Jv_{1}|^{2}
                        & = & 2(|v_{1}|^{2} + |v_{2}|^{2})
                          + 4\langle v_{1}, Jv_{2} \rangle \\
                      & = &  2|df|^{2} - 4f^{*}\omega(e_{1},e_{2}).
\eest
That gives (a), and (b) follows from  the similar computation
\best
4\langle \dbar f , K(f,\alpha)\rangle
& = &\langle v_{1} + J v_{2}  ,
                  K_{\alpha} v_{2} - J K_{\alpha} v_{1} \rangle
        + \langle v_{2} - J v_{1}  ,
                  -K_{\alpha} v_{1} - J K_{\alpha} v_{2} \rangle \\
& = & \langle v_{1}  ,   K_{\alpha} v_{2} \rangle
        - \langle v_{1}  ,   JK_{\alpha} v_{1} \rangle
        + \langle Jv_{2}  ,   K_{\alpha}\, v_{2} \rangle
        - \langle Jv_{2}  ,   JK_{\alpha} v_{1} \rangle \\
& &\hspace{.4in}  - \langle v_{2}  ,   K_{\alpha}v_{1} \rangle
        - \langle v_{2}  ,   JK_{\alpha} v_{2} \rangle
        + \langle Jv_{1}  ,   K_{\alpha} v_{1} \rangle
        + \langle Jv_{1}  ,   JK_{\alpha} v_{2} \rangle \\
& = & 4\langle  v_{1}  ,   K_{\alpha} v_{2} \rangle\\
       & = & 4f^{*}\alpha(e_{1}  ,  e_{2}).
\eest
Next fix $x\in X$ and an orthonormal  basis
$\{w^{1},w^{2},w^{3},w^{4}\}$ of $T^{*}_{x}X$ with
$w^{2}=-Jw_{1}$ and $w^{4}=-Jw^{3}$.  Then the six forms
$$ w^{1}\wedge w^{2} \pm w^{3}\wedge w^{4},\ \
         w^{1}\wedge w^{3} \pm w^{2}\wedge w^{4},\ \
         w^{1}\wedge w^{4} \pm w^{2}\wedge w^{3} $$
give an orthonormal basis of $\Lambda^{2}(T^{*}_{x}X)$,
       and two of these span the subspace of $J$ anti-invariant forms.  Hence
$$ \alpha = a(w^{1}\wedge w^{3} - w^{2}\wedge w^{4}) +
                  b(w^{1}\wedge w^{4} + w^{2}\wedge w^{3}) $$
for some $a$ and $b$, and in this basis  $K_{\alpha}$ is the matrix
\begin{equation*}
       \left(
        \begin{matrix}
            0  &  0  &  a  &  b  \\
            0  &  0  &  b  &  -a  \\
           -a  & -b  &  0  &  0  \\
           -b  & a  &  0  &  0
        \end{matrix}
       \right)
\end{equation*}
Consequently,   $K_{\alpha}^{2} = -(a^{2} + b^{2})Id =
-|\alpha|^{2}Id$.  Lastly, since $K_\alpha$ is skew-adjoint, (c) shows that
$$
|K_{J}(f,\alpha)|^{2}\, =\, -\langle \partial_{J} f\circ j,\
K_\al^2(\partial_{J} f\circ j)\rangle
=\,\,f^{*}(\,|\al|^2)\,|\partial_{J} f|^2.
$$
Equation (d) then follows from (a) because
$|df|^2=|\partial_{J} f|^2+|\dbar f|^2$.
\QED

\bigskip


\begin{cor}\label{C:Basic}
Suppose the map $f:\Sigma\to X$ is $(J,\alpha)$-holomorphic.
Then
\begin{enumerate}
\item[(a)] $|\dbar f|^{2}\ dv = f^{*}\alpha$,
\item[(b)] $\left(1-f^{*}(\,|\alpha|^{2})\,\right)\,\,|df|^{2}\ dv
                 = 2\,\left(1+f^{*}(|\alpha|^{2})\,\right)\,
                   f^{*}\omega$, \ \ and
\item[(c)] $f^{*}(|\alpha|^{2})\,|df|^{2}
                 = \left(1+f^{*}(|\alpha|^{2})\,\right)\,|\dbar f|^{2}$.
\end{enumerate}
\end{cor}

\pf Since $f$ is  $(J,\alpha)$-holomorphic,
$|\dbar f|^{2} =
        \langle \dbar f , K_{J}(f,\alpha) \rangle
        = | K_{J}(f,\alpha) |^{2}$, so (a) follows from
Proposition~\ref{P:Basic}b
while (b) and (c) follow from Proposition~\ref{P:Basic} (a) and (d).
\QED

\vskip1.5cm

       There is a second way of writing the ($J,\alpha$)-holomorphic equation
(\ref{1.defJalphaholo}).   For each $\alpha\in\Omega^{-}_{J}(X)$,
$Id + JK_{\alpha}$ is invertible
since $J K_{\alpha}$ is skew-adjoint.
Hence
\begin{equation}
\label{1.defJalpha}
J_{\alpha}=(Id + JK_{\alpha})^{-1} J(Id + JK_{\alpha})
\end{equation}
is an  almost complex structure. A map
$f:\Sigma \to X$ is $(J,\alpha)$-holomorphic
if and only if $f$ is $J_{\alpha}$-holomorphic, i.e. satisfies
\begin{equation}\label{E:CR}
\dbara f\ =\ \frac12\left(df + J_{\alpha}dfj\right)\ =\  0.
\end{equation}
      From this perspective, a solution of
the ($J,\alpha$)-holomorphic equation is a
$J_{\alpha}$ holomorphic map
with $J_{\alpha}$ lying in the family  (\ref{1.defJalpha})
parameterized by
       $\alpha\in\Omega^{-}_{J}(X)$.  In particular, we see from
(\ref{E:CR}) that the
($J,\alpha$)-holomorphic equation is elliptic.


\begin{proposition}\label{P:Basic2}
For any $\alpha\in\Omega^{-}_{J}(X)$,
the almost complex structure $J_{\alpha}$ on $X$ satisfies
\bear
\label{P:Basic2eq}
       \langle J_{\alpha}u , J_{\alpha}v \rangle
                   = \langle u,v \rangle \qquad {\mbox{and}}\qquad
       J_{\alpha} = \dfrac{1-|\alpha|^{2}}{1+|\alpha|^{2}}\ J
                              -\dfrac{2}{1+|\alpha|^{2}}\ K_{\alpha}
\eear
\end{proposition}

\pf From (\ref{1.3}), the endomorphisms
$ A_{+} = Id + JK_{\alpha} $ and $ A_{-} = Id - JK_{\alpha} $
are transposes, and $A_{+}J = JA_{-}$ and
$A_{+}K_{\alpha} = K_{\alpha}A_{-}$.
Consequently, $A_{+}^{-1}$ and $A_{-}^{-1}$ are transposes, with
$A_{-}^{-1}J = JA_{+}^{-1}$ and
$A_{-}^{-1}K_{\alpha} = K_{\alpha}A_{+}^{-1}$ and therefore
$A_{-}^{-1}A_{+} = A_{+}A_{-}^{-1}$.
Consequently,
\begin{align*}
\langle J_{\alpha}u , J_{\alpha}v \rangle
        & = \langle A_{+}^{-1}JA_{+}u , A_{+}^{-1}JA_{+}v \rangle
          = \langle JA_{-}^{-1}A_{+}u , JA_{-}^{-1}A_{+}v \rangle  \\
        & = \langle A_{-}^{-1}A_{+}u , A_{-}^{-1}A_{+}v \rangle
          = \langle u , A_{-}A_{+}^{-1}A_{-}^{-1}A_{+}v \rangle  \\
        & = \langle u , v \rangle .
\end{align*}
On the other hand, noting that $K_{\alpha}^{2} = -|\alpha|^{2}Id$, it
is easy to verify that
\begin{equation}\label{E:Inverse}
(Id + JK_{\alpha})^{-1} = \dfrac{1}{1 + |\alpha|^{2}}\,Id
                              -\dfrac{1}{1 + |\alpha|^{2}}\,JK_{\alpha}.
\end{equation}
With that, the second part of (\ref{P:Basic2eq}) follows from the
definition of $J_{\alpha}$.
\QED

\bigskip

In summary,  ($J,\alpha$)-holomorphic maps can be
regarded as solutions of  the
$J_{\alpha}$-holomorphic map equation $\dbara f=0$ for a family of
almost complex structures
parameterized by $\al$ as in (\ref{E:CR}).  We will frequently move
between these two viewpoints.

\vskip1cm


\setcounter{equation}{0}
\section{Curves and Canonical Families of $(J,\al)$ Maps}
\label{section2}
\bigskip

Given a K\"{a}hler surface $X$, we would like to use
$(J,\al)$-holomorphic curves to solve the following
problem in   enumerative geometry:

\bigskip

\noindent{\bf Enumerative Problem}\ \ Give  a $(1,1)$ homology class $A$, count
the  curves in $X$ that represent
$A$, have a specified genus $g$, and pass through the  appropriate
number of generic points.

\bigskip

We begin this section with some dimension counts which show that in
order to interpret this problem in
terms of holomorphic maps we need to consider families of maps of
dimension $p_g$.  We then show that
there is a very natural family of $(J,\al)$-holomorphic maps with exactly
that many parameters.  We conclude
the section with a theorem showing that such maps do indeed
represent holomorphic curves in $X$.

\bigskip

One can phrase the above enumerative problem in  terms of
the ${\em Severi\ variety}$\ $V_{g}(C)\subset |C|$, which
is defined  to be the closure of the set of all curves with
geometric genus $g$.  Assuming that $C-K$ is ample, it  follows from the
Riemann-Roch theorem that the dimension of the complete linear system $|C|$
is given in terms of $p_{g}=\mbox{dim}_{\,{\Bbb C}}H^{0,2}(X)$ and
$q=\mbox{dim}_{\,{\Bbb C}}H^{0,1}(X)$ by
\begin{equation*}
\mbox{dim}_{\,{\Bbb C}}|C|\ =\
\frac{C^{2}-C\cdot K}{2}\  +\  p_{g}\  -\ q
\end{equation*}
and the formal dimension of the  Severi variety is
\begin{equation}
\label{E:dim-SV}
\mbox{dim}_\cx V_{g}(C)\  =\
-K\cdot C\  +\  g - 1 \ +\ p_{g}\ -\ q.
\end{equation}
The right-hand side of (\ref{E:dim-SV}) is the `appropriate number'
of point constraints to impose; the
set of curves in $V_g(C)$ through that many generic points should be
finite, making the enumerative
problem well-defined.

\smallskip

Now let ${\cal M}_{g}(X,A)$
be the moduli space of
holomorphic maps from Riemann surfaces of genus $g$,
which represent homology class $A$.
Then its virtual dimension is given by
\begin{equation}
            \label{E:dim-MC}
\dim_{\,{\Bbb C}}{\cal M}_{g}(X,A)\ =\
-K\cdot A\  +\  g - 1.
\end{equation}

\medskip
The image of a map in  $\overline{\cal M}_{g}(X,[C])$
might be not a divisor in $|C|$, instead
it is a divisor in some other complete linear system
$|C^{\prime}|$ with $[C^{\prime}]=[C]$.
As in [BL3], we define the parameterized Severi variety
\begin{equation*}
V_{g}([C]) =
         \underset{[C^{\prime}]=[C]}{\coprod} V_{g}(C^{\prime})
\end{equation*}
Its expected dimension is now given by
\begin{equation}
            \label{E:dim-PSV}
\mbox{dim}_{\,{\Bbb C}}V_{g}([C])\ =\
               -K\cdot C \ + \ g-1\ +\ p_{g}.
\end{equation}
We still have $p_{g}$ dimensional difference between
(\ref{E:dim-PSV}) and (\ref{E:dim-MC}).
Therefore, the cut-down moduli space
by (\ref{E:dim-PSV}) many point constraints is empty
when $p_{g}\geq 1$.
This implies that the GW invariants
with (\ref{E:dim-PSV}) many point constraints are zero,
whenever $p_{g}\geq 1$.

\vskip.3in

We show that there is  a natural --- in fact obvious ---
$p_g$-dimensional family of
$(J,\al)$-holomorphic maps associated with every K\"{a}hler surface.

\begin{defn}
Given a K\"{a}hler surface
$X$, define the parameter space ${\cal H}$ by
\begin{equation}\label{D:Parameter}
{\cal H} \ =\left\{ \al+\overline{\al}\ \left|\  \al\in
H^{2,0}(X)\right\}\right.
\end{equation}
\label{2.H}
\end{defn}
Here $H^{2,0}(X)$ means the set of holomorphic $(2,0)$ forms on $X$.
Note that all forms $\al\in
H^{2,0}(X)$ are closed since $d\al=\partial \al
+\overline{\partial}\al=\partial\al$ is a $(3,0)$ form and hence vanishes
because $X$ is a complex surface.  Thus  ${\cal H}\subset \Omega^-_J(X)$ is a
$2p_{g}$-dimensional  real  vector space of closed forms.  We give it
the (real)  inner product defined
by the  $L^2$ inner product of forms:
\bear
\label{2.innerproduct}
\langle\al,\ \beta\rangle\ =\ \int_X \al\wedge\beta.
\eear

     We can use the forms $\al\in {\cal H}$ to
parameterize the right-hand side of the $(J,\al)$-holomorphic map
equation (\ref{def1.2}).

\begin{defn}
Henceforth the term `$(J,\al)$-holomorphic map' means a map
satisfying (\ref{def1.2}) for $\al$ in the
above family ${\cal H}$.
\label{2.jholo}
\end{defn}

\begin{lemma}
\label{L:ZS}
The zero divisor $Z(\alpha)$ of each nonzero
$\alpha\in {\cal H}$ represents the
canonical class.
\end{lemma}

\pf
Write $\alpha =
              \beta + \overline{\beta}$ with $\beta\in H^{2,0}(X)$.
Since $\beta$ is a section of the
canonical bundle, this means that
$Z(\alpha)=Z(\beta)$ represents the canonical divisor
with appropriate multiplicities.
\QED

\bigskip
Next, using this $2p_{g}$ dimensional parameter space $\PH$,
we define the family moduli space
\begin{equation*}
\overline{\cal M}_{g}^{\PH}(X,[C])\ =\
\{\ (f,\alpha)\ |\ \overline\partial_{J_{\alpha}}f=0,\
                       [\mbox{Im}\,f]=[C],\
                       \mbox{and}\ \alpha\in \PH\ \}
\end{equation*}
Since we just parameterize the $\overline\partial$-operator
by $2p_{g}$ dimensional parameter space, the formal dimension
of the family moduli space is given by
\begin{equation*}
\mbox{(\,Formal\,)\ }
\mbox{dim}_{\,{\Bbb C}}\overline{\cal M}_{g}^{\PH}(X,[C])\ =\
-K\cdot C\ + \ g-1\ +\ p_{g}
\end{equation*}

\medskip
On the other hand,
we define {\em a component of the canonical class}
to be a homology class of a component
of some canonical divisor.

\begin{theorem}\label{L:GT}
If $f$ is a $(J,\al)$-holomorphic map which represents a class $A\in
H^{1,1}(X)$.   Then $f$ is, in fact, holomorphic.
Furthermore, if $A$ is not a linear combination of
components of the canonical class, then $\alpha = 0$.
\end{theorem}

\pf
Since $\alpha\in H^{2,0}(X)\oplus H^{2,0}(X)$ is  closed   and  $A\in
H^{1,1}(X)$, it follows from Corollary~\ref{C:Basic}a that
$$\int_{\Sigma}|\overline{\partial}_{J}f|^{2}
         = \alpha(A)=0. $$
Thus $f$ is holomorphic, that is, $\dbar f \equiv 0$.
Consequently, $ |\alpha|^{2}|df|^{2} \equiv 0 $ by
Corollary~\ref{C:Basic}c.
Since $df$ has at most finitely many zeros,
we can conclude that $\al=0$ along the image of $f$.
Hence $\alpha= 0$, otherwise it contradicts to
the assumption on $A$ by Lemma~\ref{L:ZS}.
\QED

\vskip1cm


\setcounter{equation}{0}
\section{Family GW-Invariants}
\label{section3}
\bigskip

Let $X$ be a complex surface with a K\"{a}hler structure
$(\omega,J,g)$.  In this section we will define the Family
Gromov-Witten Invariants
associated to $(X,J)$ and the parameter space ${\cal H}$ of
(\ref{D:Parameter}).  We also state some
properties of these invariants.

    Our approach is the same  analytic arguments as that of Li and Tian
\cite{LT} to show that  the moduli space of $(J,\alpha)$-holomorphic
maps carries a virtual
fundamental class whenever it is compact.  While compactness is
automatic for  the usual
Gromov-Witten invariants, it must be verified case-by-case for the family GW
invariants (see Example~\ref{Ex:CPT}).  Thus compactness appears as a hypothesis
in the results of this section.

\bigskip

First, we  recall the notion of $C^\ell$ stable maps as defined in
\cite{LT}. Fix an integer
$l\geq 0$ and consider pairs $(f;\Sigma, x_{1},\cdots,x_{k})$ consisting
of
\begin{enumerate}
     \item a connected nodal curve $\Sigma =
\bigcup\limits_{i=1}^{m}\Sigma_{i}$ of arithmetic genus $g$
          with   distinct smooth marked points $x_1,\cdots,x_{k}$, and
     \item a continuous map $f:\Si\to X$ so that each restriction
           $f_i=f_{|_{\Sigma_{i}}}$ lifts to a $C^{l}$-map
           from the normalization $\tilde{\Sigma_{i}}$
           of $\Sigma$ into X.
\end{enumerate}

\begin{defn}
A {\em stable $C^{l}$ map} of genus $g$
with $k$ marked points is a pair $(f;\Sigma, x_{1},\cdots,x_{k})$
as above which  satisfies the  stability condition:

\medskip
\noindent$\bullet$\ \  If the homology class  $[f_i]\in
H_{2}(X,{\Bbb Q})$ is trivial,
           then the number of marked points in $\Sigma_i$ plus
           the arithmetic genus of $\Sigma_{i}$ is at  least three.
\medskip
\end{defn}

Two stable maps $(f,\Sigma; x_{1},\cdots,x_{k})$ and
$(f^{\prime},\Sigma^{\prime}; x_{1}^{\prime},\cdots,x_{k}^{\prime})$
are equivalent if there is a biholomorphic map $\sigma:\Sigma\mapsto
\Sigma^{\prime}$ such that
$\sigma(x_{i})=x_{i}^{\prime}$ for $1\leq i\leq k$ and
$f^{\prime}=f\circ\sigma$.
We denote by
$$
\F
$$
the space of all equivalence classes $[ f; \Sigma,  x_{1},\cdots,x_{k}]$
of $C^{l}$-stable maps of genus $g$ with $k$ marked points
and with total homology class $A$.
The topology of
$\F$
is defined by sequential convergence as in section 2   of \cite{LT}.
There are two continuous maps
from $\overline{{\cal F}}^{l}$.  First, there is an  evaluation map
\begin{equation}
ev :  \F \to X^{k}
\label{3.defev}
\end{equation}
which records the images of the marked points.
Second, for  $2g + k \geq 3$, collapsing the unstable components of
the domain gives a stabilization
map
\begin{equation}
st : \F \to
\overline{{\cal M}}_{g,k}
\label{3.defst}
\end{equation}
to the compactified  Deligne-Mumford space of genus $g$ curves with
$k$ marked points.
For $2g + k < 3$ we define $\overline{{\cal M}}_{g,k}$  to be the
topological space of consisting of a single  point and define
(\ref{3.defst}) to be the map to
that point.

\smallskip

We next  construct
a `generalized  bundle' $E$ over  $\F\times {\cal H}$, again
following  \cite{LT}.  Recall that each
$\alpha\in {\cal H}$ defines an almost complex
structure $J_{\alpha}$ on $X$ by (\ref{1.defJalpha}).
Denote by ${\rm{Reg}}(\Sigma)$ the set of all
smooth points of $\Sigma$.
For each $([ f; \Sigma,  x_{1},\cdots,x_{k}],\alpha)$,
define
$$ \Lambda_{j_{\Sigma}J_{\alpha}}(f^{*}TX) $$
to be the set of all continuous sections $\nu$ of
${\rm Hom}(T{\rm Reg}(\Sigma),f^{*}TX)$
with
$\nu \circ j_{\Sigma} = -J_{\alpha}\circ \nu$
    such that
$\nu$  extends continuously
across the nodes of $\Sigma$.
We take $E$ to be  the bundle whose fiber over $([ f,\Sigma;
x_{1},\cdots,x_{k}],\alpha)$
is $\Lambda_{j_{\Sigma}J_{\alpha}}(f^{*}TX)$
and give  $E$  the continuous topology as in section 2 of \cite{LT}.
We then define a section
$\Phi : \F
    \times {\cal H} \to E$ by
\begin{equation}\label{D:CRE}
\Phi([ f,\Sigma; x_{1},\cdots,x_{k}],\alpha)
     = df + J_{\alpha}dfj_{\Sigma}.
\end{equation}

\medskip

The right-hand side of (\ref{D:CRE}) vanishes for
$J_\alpha$-holomorphic maps.  Thus $\Phi^{-1}(0)$
is the moduli space of $(J,\alpha)$-holomorphic maps. The following
is a  family version of
Proposition 2.2 in \cite{LT}.


\begin{proposition}
\label{P:VMC}
Suppose that the set $\Phi^{-1}(0)$ is compact.
Then the section $\Phi$ gives rise to a
generalized Fredholm orbifold bundle with
a natural orientation and with index
\bear
\label{2.index}
r\ =\  2\,c_{1}(X)\,[A]\  +\  2\,(g-1)\  +\  2\,k\  +\
{\rm dim }\ \H.
\eear
\end{proposition}

\medskip

By Theorem 1.2 of \cite{LT}, the bundle $E$ has a
rational homology ``Euler class'' in $\F\times {\cal H}$; in fact,
since  ${\cal H}$ is contractible this Euler class lies in
$H_{r}(\F;{\Bbb Q})$ where
$r$ is the index (\ref{2.index}).
We call this class the {\em virtual fundamental cycle}
of the moduli space of
family holomorphic maps parameterized by ${\cal H}$ and denote it by

  \begin{equation}
    \label{3.defVFC}
    [\,\CM^{J,\H}_{g,k}(X,A)\,]^{\rm vir}.
  \end{equation}
In particular,
  \begin{equation}
    \label{3.dimformula}
    {\rm dim }\  [\,\CM^{J,\H}_{g,k}(X,A)\,]^{\rm vir}
    \ =\ 2\,c_{1}\,(X)[A]\  +\  2\,(g-1)\  +\  2\,k\  +\  2\,p_g.
  \end{equation}

\medskip

\medskip

The next issue is whether  the  virtual fundamental cycle is
independent of the  K\"{a}hler
structure on $X$. The next proposition is analogous to
the Proposition 2.3 in \cite{LT}.  It shows that  the  virtual
fundamental cycle depends only on certain
deformation class of the   K\"{a}hler structure.


\begin{proposition}
\label{P:Cobo}
Let $(\omega_{t},J_{t},g_{t})$, $0\leq t \leq 1$, be a continuous family of
K\"{a}hler structures on $X$.  Let
${\cal H}_{t}$ be the corresponding  continuous family of finite
subspaces defined by
(\ref{D:Parameter}) and let  $\Phi_{t}$ be the corresponding family
of sections of $E_{t}$ over
$\F\times {\cal H}_{t}$.
If $\Phi^{-1}_{t}(0)$ is compact for all $0\leq t \leq 1$, then
$$
[\,\CM^{J_{0},{\cal H}_{0}}_{g,k}(X,A)\,]^{\rm vir}
\ =\  [\,\CM^{J_{1},{\cal H}_{1}}_{g,k}(X,A)\,]^{\rm vir}.
$$
\end{proposition}

\medskip

The family GW invariants can now be defined by pulling back
cohomology classes by the evaluation and
    stabilization maps and integrating over the virtual fundamental
cycle.  That of course requires that
the virtual fundamental cycle exists, so we must assume that we
are in a situation where
$\Phi^{-1}_{t}(0)$ is compact.

\begin{defn}
\label{defn3.4}
Whenever  the virtual fundamental cycle
$[{\cal M}^{J,{\cal H}}_{g,k}(X,A)]^{\rm vir}$ exists,
we define the family GW invariants of $(X,J)$ to be the map
$$ GW^{J,{\cal H}}_{g,k}(X,A) :
      \left[H^{*}(X;{\Bbb Q})\right]^{k} \times H^{*}(\overline{{\cal
M}}_{g,k};{\Bbb Q})
      \mapsto {\Bbb Q} $$
defined on
$\alpha\in [\,H^{*}(X;{\Bbb Q})\,]^{k}$
and  $\beta\in
H^{*}(\overline{{\cal M}}_{g,k};{\Bbb Q})$ by
$$ GW^{J,{\cal H}}_{g,k}(X,A)(\beta;\alpha)\ =\
      [\,\CM^{J,{\cal H}}_{g,k}(X,A)\,]^{\rm vir} \cap (\,
      st^{*}\beta\cup ev^{*}\alpha\,).
$$
\end{defn}

\bigskip

The condition that $\Phi^{-1}(0)$ is compact must be checked ``by
hand''.  In general,
$\Phi^{-1}(0)$ is compact for some choices of $A$, but not for others.

\medskip

\begin{ex}\rm\label{Ex:CPT}
 Let $(X,J)$ be a  K\"{a}hler surface with $p_g>1$.  Then there is a
    non-zero element $\beta\in H^{2,0}$ whose zero set $Z(\beta)$ is
non-empty,  represents the
canonical class $K$, and whose irreducible components can be
parameterized by  holomorphic maps.
Fix a parameterization $f:\Sigma\to X$ of one such component; this
represents a non-zero class $A\in
H_2(X,\Z)$.  Then  $\alpha=\beta+\overline{\beta}$ lies
in the space ${\cal H}$ of (\ref{2.H}) and
$\Phi(f,\lambda\alpha)=0$ for all real
$\lambda$.  Thus on any  K\"{a}hler surface with $p_g>1$, the set
$\Phi^{-1}(0)$ is not compact for an component of the canonical class $A$.
\smallskip

    On the other hand,  in the next section we will give examples of
classes $A$ in  K\"{a}hler
surfaces with $p_g>1$  for which $\Phi^{-1}(0)$ {\em is} compact.
\label{noncompactremark}
\end{ex}

\medskip

Recall that the energy of a map $f:(\Sigma,j)\to X$ is defined by
$$
  E(f)\ =\ \frac{1}{2}\int_{\Sigma}\,|df|^{2}\,d\mu
$$
where $\mu$ is the metric in the conformal class $j$.

\medskip
\begin{theorem}
\label{theorem3.6}
If there is a constant $C$, depending only on $(X,\w,J,g)$  such that
$E(f)+||\alpha||<C$ for all
$(J,\alpha)$-holomorphic maps into $(X,J)$, then $\Phi^{-1}(0)$  is
compact and hence the
family  GW invariants are well-defined.
\end{theorem}

\pf Consider a sequence $(f_n, \alpha_n)$ of $J_\alpha$-holomorphic maps.
       The uniform bound on $||\alpha_n||$  implies that  the $J_\alpha$
lie in a compact
family. Since $E(f_n)<C$ the proof of Gromov's Compactness Theorem
(see \cite{pw} and \cite{is})
shows that $\{(f_n,\alpha_n)\}$ has a convergent subsequence.  Consequently,
$\Phi^{-1}(0)$ is compact as in the  hypothesis of Proposition
\ref{P:VMC}.   That means that the
    virtual fundamental cycle (\ref{3.defVFC}) is well-defined.  The
family  GW invariants are then
given by Definition \ref{defn3.4}.
\QED

\bigskip


\smallskip
We conclude this section
by an important property of the
family GW invariants.
This property  generalizes the composition law of  ordinary
Gromov-Witten invariants.
For that we consider maps from a  domain $\Sigma$ with node  $p$ and
relate them to maps whose domain
is the normalization of $\Sigma$ at $p$.  When the node is separating
the genus and the number of
marked points decompose as $g=g_{1}+g_{2}$ and $k=k_{1}+k_{2}$ and
is a natural map
\begin{equation}\label{gluing1}
\sigma : \overline{\cal M}_{g_{1},k_{1}+1}
               \times \overline{\cal M}_{g_{2},k_{2}+1}
               \mapsto \overline{\cal M}_{g,k}
\end{equation}
defined by gluing $(k_{1}+1)$-th marked point of the first component
to the first marked point of the second component.
We denote by $PD(\sigma)$ the Poincar\'{e} dual of the image of this
map $\sigma$.

\smallskip
We denote by
$E_{1}\oplus \tilde{E_{2}}$ ( resp.  $\tilde{E_{1}}\oplus E_{2}$ )
the generalized bundle over
  \begin{equation*}
    \overline{{\cal F}}_{g_{1},k_{1}+1}(X,A_{1}) \times
    \overline{{\cal F}}_{g_{2},k_{2}+1}(X,A_{2}) \times
    {\cal H} \times [0,1]
  \end{equation*}
whose fiber over
$(\,[ f_{1},\Sigma_{1}; \{x_{i}\}],
         [ f_{2},\Sigma_{2}; \{y_{j}\}], \alpha ,t\,)$
is $\Lambda^{0,1}_{j_{\Sigma_{1}}J_{\alpha }}\oplus
        \Lambda^{0,1}_{j_{\Sigma_{2}}J_{t\alpha }}$
   ( resp.
    $\Lambda^{0,1}_{j_{\Sigma_{1}}J_{t\alpha }}\oplus
        \Lambda^{0,1}_{j_{\Sigma_{2}}J_{\alpha }}$ ).
The formula
\begin{align*}\label{section-com}
    \Psi_{1}
         (\,[ f_{1},\Sigma_{1}; \{x_{i}\}],
         [ f_{2},\Sigma_{2}; \{y_{j}\}],  \alpha \,)
 \ &=\  (\, df_{1} + J_{\alpha}df_{1}j_{\Sigma_{1}},
        df_{2} + J_{t\alpha}df_{2}j_{\Sigma_{2}}\, )\\
( \mbox{\ resp.\ \ }
     \Psi_{2}
         (\,[ f_{1},\Sigma_{1}; \{x_{i}\}],
         [ f_{2},\Sigma_{2}; \{y_{j}\}],  \alpha \,)
 \ &=\  (\, df_{1} + J_{\alpha}df_{1}j_{\Sigma_{1}},
        df_{2} + J_{t\alpha}df_{2}j_{\Sigma_{2}}\, )\ )
\end{align*}
then defines a section of $E_{1}\oplus \tilde{E_{2}}$
( resp. $\tilde{E_{1}}\oplus E_{2}$ ).

\bigskip
On the other hand, for non-separating nodes there is another natural map
\begin{equation}\label{gluing2}
 \theta : \overline{\cal M}_{g-1,k+2}
              \mapsto  \overline{\cal M}_{g,k}
\end{equation}
    defined by gluing the last two marked points.
We also write $PD(\theta)$ for the Poincar\'{e} dual
of the image of  $\theta$.  The composition law is then the following
    two formulas.

\medskip
\begin{proposition}[Composition Law]
\label{P:CL}
Let $\{ H_{\gamma}\}$ be any basis of $H^{*}(X;Z)$ and
$\{ H^{\gamma}\}$ be its dual basis and suppose that
$GW^{J,\H}_{g,k}(X,A)$ is defined.
\begin{enumerate}
\item [(a)]
  For any decomposition $A=A_{1}+A_{2}$\
  if either \  $\Psi^{-1}_{1}(0)$\  or \
  $\Psi^{-1}_{2}(0)$\  is compact, then
      \begin{align*}
        &GW^{J,{\cal H}}_{g,k}(X,A)(PD(\sigma);\alpha)\
        =
       \sum_{A=A_{1}+A_{2}}I_{A_{1},A_{2}}(\,\alpha\,)
        \mbox{\ \ \ \ where\ \ }\\
    &I_{A_{1},A_{2}}(\,\alpha\,)=
       \sum_{\gamma}GW^{J,{\cal H}}_{g_{1},k_{1}+1}(X,A_{1})
         ( \alpha_{1},H_{\gamma})
         GW_{g_{2},k_{2}+1}(X,A_{2})(H^{\gamma},\alpha_{2})
         \mbox{\ if\ }
          \Psi_{1}^{-1}(0) \mbox{\ is\ compact,}   \\
    &I_{A_{1},A_{2}}(\,\alpha\,)=
         \sum_{\gamma}GW_{g_{1},k_{1}+1}(X,A_{1})
         ( \alpha_{1},H_{\gamma})
         GW^{J,\H}_{g_{2},k_{2}+1}(X,A_{2})(H^{\gamma},\alpha_{2})
         \mbox{\ if\ }
          \Psi_{2}^{-1}(0)
         \mbox{\ is\ compact,}
    \end{align*}
  $GW$ denotes the ordinary GW invariant,
  and $\alpha=\alpha_{1}\otimes\alpha_{2}$ in
  $[\,H^{*}(X)\,]^{k_{1}}\otimes [\,H^{*}(X)\,]^{k_{2}}$
\item[(b)]
  $GW^{J,{\cal H}}_{g,k}(X,A)(PD(\theta);\alpha) =
  \sum_{\gamma}GW^{J,{\cal H}}_{g-1,k+2}(X,A)
  (\alpha,H_{\gamma},H^{\gamma}).$
\end{enumerate}
\end{proposition}

\medskip
\noindent{\bf Sketch of proof\ \ }

\smallskip
\noindent
{\bf (a)\ \ }
We denote by
$E_{1}\oplus {E_{2}}$
the generalized bundle over
  \begin{equation*}
    \overline{{\cal F}}_{g_{1},k_{1}+1}(X,A_{1}) \times
    \overline{{\cal F}}_{g_{2},k_{2}+1}(X,A_{2}) \times
    {\cal H}
  \end{equation*}
whose fiber over
$(\,[ f_{1},\Sigma_{1}; \{x_{i}\}],
         [ f_{2},\Sigma_{2}; \{y_{j}\}], \alpha \,)$
is $\Lambda^{0,1}_{j_{\Sigma_{1}}J_{\alpha }}\oplus
        \Lambda^{0,1}_{j_{\Sigma_{2}}J_{\alpha }}$.
The formula
\begin{equation*}\label{section-com}
    \Psi_{1}
         (\,[ f_{1},\Sigma_{1}; \{x_{i}\}],
         [ f_{2},\Sigma_{2}; \{y_{j}\}],  \alpha \,)
 \ =\  (\, df_{1} + J_{\alpha}df_{1}j_{\Sigma_{1}},
        df_{2} + J_{\alpha}df_{2}j_{\Sigma_{2}}\, )
\end{equation*}
also defines a section of $E_{1}\oplus {E_{2}}$.
Similarly to Proposition~\ref{P:VMC},
this bundle is a generalized Fredholm orbifold bundle.
We denote its virtual fundamental cycle by $J_{A_{1},A_{2}}$.

\smallskip
On the other hand,
there is a natural map
\begin{equation*}
p : \overline{{\cal F}}_{g_{1},k_{1}+1}(X,A_{1})\times
      \overline{{\cal F}}_{g_{2},k_{2}+1}(X,A_{1})\times \PH
  \to X\times X
\end{equation*}
defined by
$(\,\left[f_{1},\Sigma_{1};\{x_{i}\}\right],
  \left[f_{2},\Sigma_{2};\{y_{j}\}\right],\alpha\,)
  \to (\,f_{1}(x_{k_{1}+1}),f_{2}(y_{1})\,)$.
Thus we have a homology class
\begin{equation*}
  \sum\,\sum_{\gamma}
  J_{A_{1},A_{2}}\
  \cap\ (\,ev_{k_{1}+1}^{*}H_{\gamma}\cup ev_{1}^{*}H^{\gamma}\,)
  \mbox{\ \ in\ \ }H_{*}(\,\cup\,p^{-1}(\triangle);{\Bbb Q}\,)
\end{equation*}
where the sum and the union are over all decompositions of $A$ and
$\triangle$ is the diagonal in $X\times X$.
We set
\begin{equation}\label{itm}
  I_{A_{1},A_{2}}\ =\
  \sum_{\gamma}
  J_{A_{1},A_{2}}\
  \cap\ (\,ev_{k_{1}+1}^{*}H_{\gamma}\cup ev_{1}^{*}H^{\gamma}\,)
\end{equation}

\smallskip
There is also a surjective map
  \begin{equation*}
   \pi : \cup\,p^{-1}(\triangle)
   \to st^{-1}(\mbox{Im}\,\sigma)
  \end{equation*}
obtained by
identifying $x_{k_{1}+1}$ and $y_{1}$,
where
$st$ is the stabilization map of $\F\times\PH$.
Note that we can consider
\begin{equation*}
[\,\CM_{g,k}^{J,\PH}(X,A)\,]^{\mbox{vir}}\cap \mbox{PD}(\sigma)
\end{equation*}
as a homology class in
$H_{*}(st^{-1}(\mbox{Im}\,\sigma);{\Bbb Q})$.
Similarly to the case of GW invariants, we have
\begin{equation}\label{fcl-2}
[\,\CM_{g,k}^{J,\PH}(X,A)\,]^{\mbox{vir}}\ \cap\
 \mbox{PD}(\sigma)
=\ \pi_{*}(\ \sum\,I_{A_{1},A_{2}}\ ).
\end{equation}

\medskip
Now, suppose $\Psi_{1}^{-1}(0)$ is compact.
For each $0\leq t\leq 1$, consider
the restriction of $E_{1}\times\tilde{E_{2}}$  and its section
$\Psi_{1}$ to
  \begin{equation*}\label{restriction}
    \overline{{\cal F}}_{g_{1},k_{1}+1}(X,A_{1}) \times
    \overline{{\cal F}}_{g_{2},k_{2}+1}(X,A_{2}) \times
    {\cal H} \times \{t\}.
  \end{equation*}
This is also a generalized Fredholm orbifold bundle with
a virtual fundamental cycle, denoted by
\begin{equation*}
[\,\CM_{t}\,]^{\mbox{vir}} =
[\,\CM_{(g_{1},k_{1}+1),
      (g_{2},k_{2}+1)}^{J,\PH}(X,A_{1},A_{2},t)\,]^{\mbox{vir}}.
\end{equation*}
Note that by definition $[\,\CM_{1}\,]^{\mbox{vir}} = J_{A_{1},A_{2}}$
and
\begin{equation*}
   [\,\CM_{0}\,]^{\mbox{vir}}\  =\
[\,\CM_{g_{1},k_{1}+1}^{J,\PH}(X,A_{1})\,]^{\mbox{vir}} \otimes
[\,\CM_{g_{2},k_{2}+1}(X,A_{2})\,]^{\mbox{vir}}.
\end{equation*}
where $[\,\CM_{g_{2},k_{2}+1}(X,A_{2})\,]^{\mbox{vir}}$
is a class which defines GW invariants.
By the same argument as in Proposition~\ref{P:Cobo}, we finally have
\begin{equation*}\label{fcl-1}
J_{A_{1},A_{2}}\ =\
[\,\CM_{1}\,]^{\mbox{vir}}\  =\  [\,\CM_{0}\,]^{\mbox{vir}}\  =\
[\,\CM_{g_{1},k_{1}+1}^{J,\PH}(X,A_{1})\,]^{\mbox{vir}} \otimes
[\,\CM_{g_{2},k_{2}+1}(X,A_{2})\,]^{\mbox{vir}}
\end{equation*}
as homology classes in
$H_{*}(\,\overline{{\cal F}}_{g_{1},k_{1}+1}(X,A_{1})\times
      \overline{{\cal F}}_{g_{2},k_{2}+1}(X,A_{1});{\Bbb Q} \,).$
Together with (\ref{itm}) and (\ref{fcl-2}) this implies the first
composition law.

\medskip
\noindent
{\bf (b)\ \ }
Similarly as above, we have an evaluation map of last two marked points
\begin{align*}
p:\overline{{\cal F}}_{g-1,k+2}(X,A)\times \PH &\to X\times X  \\
  \left(\,\left[f,\Sigma;\{x_{i}\}\right],\alpha\,\right) &\to
  (\,f(x_{k+1}),f(x_{k+2})\,).
\end{align*}
There is also a surjective map
$\pi : p^{-1}(\triangle)\to st^{-1}(\mbox{Im}\,\theta)$.
It follows that
\begin{equation*}
[\,\CM_{g,k}^{J,\PH}(X,A)\,]^{\mbox{vir}}\cap PD(\theta)\  =\
\pi_{*}\left(\,\sum\,[\,\CM_{g-1,k+2}^{J,\PH}(X,A)\,]^{\mbox{vir}}\cap
             \left(\,ev_{k+1}^{*}H_{\gamma}\wedge
                     ev_{k+2}^{*}H^{\gamma}\,\right)\,\right)
\end{equation*}
which implies the second Composition Law.
\QED

\bigskip

    That completes our overview of the family GW invariants.  We next look
at some examples, namely the various types of minimal  K\"{a}hler
surfaces. There we can use the
specific geometry of the space to verify that the
     moduli space is compact and hence the family GW invariants are
well-defined.

\vskip 1cm


\setcounter{equation}{0}
\section{ K\"{a}hler surfaces with $p_{g}\geq 1$}
\label{section4}
\bigskip


In this section we will focus on the family GW-invariants
for minimal K\"{a}hler surfaces $X$ with $p_{g}\geq 1$. The
Enriques-Kodaira Classification \cite{bpv}
separates such surfaces into the following three types.
\begin{enumerate}
\item  $X$ is K3 or Abelian surface with  canonical class $K=0$.
          In this case, $p_{g}=1$.
\item  $X$ is an elliptic surface \  $\pi:X\to  C$ with
          Kodaira dimension 1.  If
          the multiple fibers $B_{i}$ have multiplicity $m_{i}$, then
          a canonical divisor is
            \begin{equation}\label{E:CD-E}
            K=\pi^{*}D +\sum(m_{i}-1)B_{i}\ \ \ \ \
            {\rm where\ \ \ } {\rm deg\ }D=2g(C)-2+\chi({\mathcal{O}}_{X})
            \end{equation}
\item $X$ is a surface of general type with $K^{2}>0$.
\end{enumerate}
We will examine these cases one at a time.  For each we will show
that the family invariants $GW^{J,{\cal H}}_{g,k}(X,A)$
are well-defined.  By Theorem \ref{theorem3.6} the  key issue is   bounding
    the energy $E(f)$ and the
pointwise norm $|\al|$   uniformly for all
$(J,\al)$-holomorphic maps into $X$.


\vskip 1cm
\noindent
{\bf  K3 and Abelian Surfaces }

\bigskip
\noindent
Let $(X,J)$ be a K3 or Abelian surface.  Since the canonical class is
trivial, Yau's proof of the Calabi
conjecture implies that $(X,J)$ has a K\"{a}hler structure
$(\omega,J,g)$ whose metric $g$ is Ricci flat.
For such a structure all holomorphic $(0,2)$ forms are parallel, and
hence have pointwise constant
norm (see \cite{b}).  Thus  $\PH\cong \cx$  consists of closed
forms $\al$ with
$|\al|$ constant.  Furthermore, the structure is also  hyperk\"{a}hler,
meaning that there is a three-dimensional space of K\"{a}hler structures
which is isomorphic as an algebra to the imaginary quaternions.   The unit
two-sphere in that space is the so-called {\em Twistor Family} of
complex structures.

\smallskip
Consider the set
${\cal T}_{0} = \{\, J_{\alpha} \ |\  \alpha\in {\cal H} \}$.
Since $\alpha$ has no zeros, equation (\ref{P:Basic2eq})  shows that
$J_{\alpha}\rightarrow -J$
uniformly as $|\alpha|\rightarrow \infty$.
We can therefore compactify ${\cal T}_{0}$ to ${\cal T}\cong \P^1$
by adding $-J$ at infinity.

\begin{proposition}\label{Twistor}
${\cal T}$ is the Twistor Family induced from the hyperk\"{a}hler
metric $g$.
\end{proposition}
\pf
Let $\alpha\in \PH$ with $|\alpha|=1$.
It then follows from Proposition~\ref{P:Basic2} that
   $J_{\alpha}=-K_{\alpha}$ and $(\alpha,J_{\alpha},g)$ is a
K\"{a}hler structure on $X$. On the other hand, we define
$\alpha^{\prime}$  by
$\alpha^{\prime}(u,v)=\alpha(u,Jv)$.
Then $|\alpha^{\prime}|=1$ and $\alpha^{\prime}\in\PH$ since
$\beta^{\prime}$ is  holomorphic for each holomorphic 2-form $\beta$.
Moreover, by definition we have
$$J_{\alpha^{\prime}}=-K_{\alpha^{\prime}}=-JK_{\alpha}=JJ_{\alpha}.$$
Since $(\alpha^{\prime},J_{\alpha^{\prime}},g)$ is also K\"{a}hler and
$JJ_{\alpha}J_{\alpha^{\prime}}=-Id$, \  the K\"{a}hler structures
$\{J,J_\alpha,J_{\alpha'}\}$  multiply as
unit imaginary quaternions.  It follows that
${\cal T}$ is the Twistor Family induced from the hyperkähler
metric $g$.
\QED

\medskip
\begin{lemma}
\label{L:K3}
Let $A$ be a nontrivial homology class with $\omega(A)\geq 0$.
Then there exits
a constant  $C_A$ such that every $(J,\al)$-holomorphic map
$f:C\to X$ representing $A$
with $\alpha\in {\cal H}$ satisfies
\begin{equation*}
E(f) = \frac{1}{2}\int_{\Sigma}|df|^{2} < \omega(A) + C_A
\qquad\mbox{and}\qquad
|\alpha| \leq 1.
\end{equation*}
\end{lemma}

\pf  Since $|\al|$ is a constant, we can integrate
Corollary~\ref{C:Basic}b to conclude that
$|\alpha|\leq 1$. Let $C_A$ be an upper bound for the  function
$\al\mapsto |\al(A)|$ on the set of
$\al\in{\cal H}$ with $|\al|\leq 1$. Because $\al$ is closed,
Proposition~\ref{P:Basic}a and
Corollary~\ref{C:Basic}a imply that
\begin{equation*}
        E(f) =
             \dfrac{1}{2}\int\limits_{C}|df|^{2}
             = \int\limits_{\Sigma}f^{*}(\omega + \alpha)
             = \omega(A) + \alpha(A) \leq \omega(A) + C_A. \mbox{\QED}
\end{equation*}

\medskip

\begin{theorem}\label{C:K3}
Let $(X,J)$ be a K3 or Abelian surface. For each non-trivial $A\in
H_2(X,\Z)$, the invariants
$GW^{J,{\cal H}}_{g,k}(X,A)$ are well-defined and independent of
     $J$.  Furthermore,  if $A=mB$ and $A^{\prime}=mB^{\prime}$
           where $B$ and $B^{\prime}$ are primitive
           with the same square,
           then $$ GW^{J,\PH}_{g,k}(X,A) =
                   GW^{J,\PH}_{g,k}(X,A^{\prime}).$$
\end{theorem}

\pf
For any nontrivial homology class $A$,
we can choose a  Ricci flat  K\"{a}hler structure $(\omega,J,g)$
   such that
$\omega(A)\geq 0$ ( if $\omega(A) < 0$, then
we choose $(-\omega,-J,g)$\,).
It then follows from
Lemma~\ref{L:K3} and
Theorem~\ref{theorem3.6} that
$GW^{J,\PH}_{g,k}(X,A)$ is well-defined.

\medskip

Bryan and Leung have applied the machinery of Li and Tian to define
family GW invariants associated to the
   Twistor Family ${\cal T}$  \cite{bl1, bl2}.  Their invariants, which
we denote by
$$
\Phi^{{\cal T}}_{g,k}(X,A),
$$
are actually independent of  the Twistor Family since the moduli
space of complex structures on $X$
is connected. On the other hand,
if $A=mB$ and $A^{\prime}=mB^{\prime}$
           where $B$ and $B^{\prime}$ are primitive
           with the same square, then
there is an orientation preserving diffeomorphism
of $X$ which sends the class $B$ to the class $B^{\prime}$.
That implies that
$\Phi^{{\cal T}}_{g,k}(X,A)=$$\Phi^{{\cal T}}_{g,k}(X,A^{\prime})$.

\medskip
To complete the proof  it suffices to show that
\begin{equation}
GW^{J,{\cal H}}_{g,k}(X,A) =
          \Phi^{{\cal T}}_{g,k}(X,A).
\label{BL=GW}
\end{equation}
For that,  recall  from Theorem 1.2 of \cite{LT} that the
moduli cycle is defined from a section $s$ of a generalized Fredholm
orbifold bundle
   $E\to B$ and is represented by a cycle that lies in an arbitrarily
small neighborhood of  $s^{-1}(0)$.  Both
sides of (\ref{BL=GW}) are defined in that way using the same
Fredholm  bundle $E$ over the space of K\"{a}hler
structures.  In the first case   $B$ is
$\{J_\alpha\,|\,\alpha\in {\cal H}\}$ and $s^{-1}(0)$ is  the set
of all $(f,\alpha)$ where $f$
is a $J_\alpha$-holomorphic map, and in the second case $B={\cal T}$
is  the Twistor Family and $s^{-1}(0)$
is  the set of
$J_\alpha$-holomorphic maps for $J_\alpha\in {\cal T}$ .  By
Proposition \ref{Twistor}
$\{J_\alpha\,|\,\alpha\in {\cal H}\}$ parameterizes the Twistor
Family after adding a  point at infinity to
${\cal H}$.  But since
$\omega(A)\geq 0$, Lemma
\ref{L:K3} shows that $|\alpha| \leq 1$ for all $J_{\alpha}$ holomorphic maps
representing the homology class $A$
with $\alpha\in {\cal H}$.  Thus the moduli cycle is bounded
away from the  point at infinity, so the two  definitions of the
moduli cycle are exactly equal.  That gives
(\ref{BL=GW})
\QED


\vskip 1cm
\noindent
{\bf  Elliptic Surfaces }

\bigskip
\noindent
First, we recall the well-known facts
about minimal elliptic surfaces $X$ with Kodaira dimension 1 \cite{fm}.
\begin{enumerate}
\item $X$ is elliptic  in a unique way.
\item Every deformation equivalence is through elliptic surfaces.
\end{enumerate}
Therefore, there is a unique elliptic structure
$\pi :( X,J) \to C$.
Moreover, for the fiber class $F$ and
any homology class $A\in H_{2}(X;Z)$, the integer
\begin{equation}\label{E:Elliptic}
F\cdot A + \mbox{deg}(\pi_{*}A)
\end{equation}
is well-defined for each complex structure $J$ and
it is invariant under the deformation
of complex structure $J$.

\medskip
Let $(\omega,J,g)$ be a K\"{a}hler structure on $X$
and $\PH$ be as in (\ref{D:Parameter}).  For $\al\in\PH$, let
$\|\al\|$ denote the $L^2$ norm as in
(\ref{2.innerproduct}).


\begin{lemma}\label{L:Elliptic}
Let $A\in H_{2}(X;Z)$ such that
the integer (\ref{E:Elliptic}) is positive.
Then, there exit uniform constants $E_0$ and $N$ such that
for any $J_{\alpha}$-holomorphic map $f:\Sigma\to X$,
representing homology class $A$,
with $\alpha\in {\cal H}$,
we have
\begin{equation*}
E(f) = \frac{1}{2}\int_{\Sigma}|df|^{2} \leq E_{A}\, ,\ \ \ \ \
||\alpha|| \leq N.
\end{equation*}
\end{lemma}
\pf
It follows from (\ref{E:CD-E})
and Lemma~\ref{L:ZS} that
for any nonzero $\alpha\in{\cal H}$, the zero set of $\al$ lies in
the union of fibers $F_i$.  Let
$N(\al)$ be a (non-empty) union of $\ep$-tubular neighborhoods of the
$F_i$.  Denote by ${\cal
S}$ the unit sphere in ${\cal H}$ and set
\begin{equation*}
        m(J) = \underset{\alpha\in {\cal S}}{\mbox {min}}\
               \underset{x\in X\setminus N(\alpha)}
               {\mbox {min}}\ |\alpha|\ \ \ \mbox{and}\ \ \
        N = \frac{2}{m(J)} .
\end{equation*}

We can always choose a smooth fiber $F\subset X\setminus N(\alpha)$
such that $f$ is transversal to $F$.
Let $f^{-1}(F)=\{p_{1},\cdots,p_{n}\}$ and for each $i$ fix   a small
     holomorphic  disk  $D_{i}$ normal to $F$  at $f(p_{i})$. We can
further assume that
$f$ is  transversal to each $D_{i}$ at $f(p_{i})$.

\smallskip
Define $\mbox{sgn}(r)$ to be the sign of a real number $r$
if $r\ne 0$, and 0 if $r=0$. Denote by
$I(S,f)_{p}$ the local intersection number of
the map $f$ and a submanifold $S\hookrightarrow X$ at $f(p)$.
In terms of bases $\{e_{1},\, e_{2}=j\,e_{1} \}$ of $T_{p_{i}}\Sigma$,
$\{v_{1},\, v_{2}=j\,v_{1} \}$ of $T_{f(p_{i})}F$, and
$\{v_3,\, v_4=j\,v_3 \}$ of $T_{f(p_{i})}D_{i}$ we have
\begin{align}
&I(F,f)_{p_{i}}  = \mbox{sgn}\left(
          (v^{1}\wedge v^{2}\wedge v^{3}\wedge v^{4})
           (v_{1},v_{2},f_{*}e_{1},f_{*}e_{2})\right)
         = \mbox{sgn} \left(
          ( v^{3}\wedge v^{4})(f_{*}e_{1},f_{*}e_{2})\right),
           \notag \\
&I(D_{i},f)_{p_{i}} = \mbox{sgn}\left(
          (v^{1}\wedge v^{2}\wedge v^{3}\wedge v^{4})
           (f_{*}e_{1},f_{*}e_{2},v_{3},v_{4})\right)
         = \mbox{sgn}\left(
         (v^{1}\wedge v^{2})(f_{*}e_{1},f_{*}e_{2})\right).
           \notag
\end{align}
Comparing with  $\mbox{sgn}\,f^{*}\omega(e_{1},e_{2}) = \mbox{sgn}\left(
          (v^{1}\wedge v^{2})(f_{*}e_{1},f_{*}e_{2}) +
          (v^{3}\wedge v^{4})(f_{*}e_{1},f_{*}e_{2}) \right)$ shows that
\bear
     I(F,f)_{p_{i}} + I(D_{i},f)_{p_{i}} \ = \
          \mbox{sgn}\ (f^{*}\omega)(e_{1},e_{2}).
     \label{E:LIN}
\eear

\smallskip
Now suppose $m(J)||\alpha|| \geq 2$.
Then $|\al|\geq 2$ along each $F_i$, so by (\ref{E:LIN}) and
Corollary~\ref{C:Basic}b
\begin{equation*}
\sum_{i}\left( I(f,F)_{p_{i}} +
                  I(f,D_{i})_{p_{i}} \right) < 0.
\end{equation*}
This contradicts to our assumption
$A\cdot f + \mbox{deg}(\pi_{*}A) > 0$
since by definition $\sum_{i}I(f,F)_{p_{i}} = A\cdot f$
and $\sum_{i}I(f,D_{i})_{p_{i}} = \mbox{deg}(\pi_{*}A).$
Therefore $||\alpha|| < N$ with $N$ as above.
The energy bound follows exactly same arguments as in the proof of
Lemma~\ref{L:K3}.
\QED

\medskip

\begin{proposition}\label{C:Elliptic}
For any homology class $A$
with (\ref{E:Elliptic}) positive,
the invariants
$GW^{J,{\cal H}}_{g,k}(X,A)$ are well-defined and depend only on the
     deformation class of $(X,J)$.
\end{proposition}
\pf
It follows from Lemma~\ref{L:Elliptic} and
Theorem~\ref{theorem3.6} that the invariants
$GW^{J,{\cal H}}_{g,k}(X,A)$ are well-defined.
On the other hand, (\ref{E:Elliptic}) is
invariant under the deformation of $J$.
Therefore, applying Proposition~\ref{P:Cobo},
we can conclude that
the invariants only depends on
the deformation equivalence class of $J$.
\QED


\vskip 1cm
\noindent
{\bf  Surfaces of General Type }

\bigskip
\noindent
Let $(X,J)$ be a minimal surface of general type.

\begin{proposition}\label{P:gt}
If $A$ is of type (1,1) and is not a linear combination
of components of the canonical class, then
we can define the invariant
$GW^{J,{\cal H}}_{g,k}(X,A)$.
They are invariant under the deformations
of complex structures which
preserve (1,1)-type of $A$.
\end{proposition}

\pf
Lemma~\ref{L:GT} and Theorem~\ref{theorem3.6} imply that
the invariants $GW^{{\cal H}}_{g,k}(X,A)$
are well-defined under the assumption that $A$ is type $(1,1)$.
On the other hand, Proposition~\ref{P:Cobo} also implies that
the invariants $GW^{{\cal H}}_{g,k}(X,A)$
are invariant under  deformations
of the complex structure which
preserve the $(1,1)$ type of $A$.
\QED


\vskip 1 cm
\noindent
{\bf A Linear Combination of Components of the Canonical Class}
\bigskip

\noindent There is an alternative way of defining invariants: As
in \cite{ip1,ip2}, we  define an invariant from a rational
homology class in $\CM_{g,k}\times X^{k}$ that is induced from the
family moduli space. This invariant is same as the family
invariant whenever both invariants are well-defined. In
particular, this invariant is also well-defined for a class $A$
which is a linear combination of components of canonical class and
satisfies (\ref{condition}) below.

\medskip
Let $(X,J)$ be a minimal K\"{a}hler surface and
$A$ be a class of type $(1,1)$ with
\begin{equation}\label{condition}
 -A\cdot K\  +\  g-1\  \ \geq\ 0
\end{equation}
where $K$ is the canonical class.
As in section 3, denote by
$E$ the generalized bundle
over  $\SF\times {\cal H}$ with a section $\Phi$
defined by $\Phi(f,\al)=df+J_{\al}\,df\,j$.

\begin{proposition}\label{P:last}
The section $\Phi$ as above
gives rise to a well-defined homology class
$[\,\CM_{g,k}^{J,\H}(X,A)\,]$ in
$H_{2r}(\,\CM_{g,k}\times X^{k};{\Bbb Q}\,)$,\
where\  \ $r=-A\cdot K +g-1 + p_{g} + k$.
\end{proposition}

\pf
It follows from Lemma~\ref{L:ZS} and Theorem~\ref{L:GT} that
\begin{equation}\label{frontier}
  \CM_{g,k}(X,A)\ \subset\
  \Phi^{-1}(0)\ \subset\ \CM_{g,k}(X,A)\
  \cup\ \left( F\times \H \right)
\end{equation}
where
$F$ is the set of all  $(f,0)$ in $\Phi^{-1}(0)$ such that
the image of $f$ lies entirely in some canonical divisor.
Since
${\rm dim}_{\Bbb C}\,|K|={\rm dim}_{\Bbb C}\,H^{2,0}(X)-1=p_{g}-1$,
under the stabilization and evaluation maps
 $st\times ev$ the image of $F\times \H$ lies in (\,real\,) dimension
less than $2\,(\,p_{g}-1+k\,)$.
Hence, the assumption (\ref{condition}) implies that
the image of $F\times \H$ lies in dimension at most $2\,r-2$.

\smallskip
Let $d$ be the dimension of \ $\CM_{g,k}\times X^{k}$.
Then there exists a neighborhood $V$ of the image of $F\times \H$
such that (i) the basis for $H_{d-2r}(\,\CM_{g,k}\times X^{k};{\Bbb Q}\,)$
are represented by cycles $D_{i}$\  for\  $1\leq i\leq m$ and (ii)
each cycle $D_{i}$ does not intersect with $V$.
Let $U$ be the preimage of $V$ under the map $st\times ev$.
It then follows from the proof of Proposition 2.2 in \cite{LT} that
$\CM_{g,k}(X,A)$ can be covered by finitely many smooth approximations
$U_{i}$\ for $1\leq i \leq n$ such that
$F \subset \cup_{l\leq i \leq n}\,U_{i} \subset U$ for some $l\leq n$.
Then using the arguments in the proof of Theorem 1.2 in \cite{LT}
we can construct a {\em cocycle} $Z$ with its boundary
lying in $U$.

As in the proof of Proposition 4.2 of \cite{km},
we now define a homology class
$[\,\CM_{g,k}^{J,\H}(X,A)\,]$ as a homology class in
$H_{2r}(\,\CM_{g,k}\times X^{k};{\Bbb Q}\,)$ determined by
the intersection numbers between $(st\times ev)(Z)$ and
the cycles $D_{i}$\  for\  $1\leq i\leq m$.
Since the image of $F$ lies in dimension at most $2\,r-2$,
this homology class is well-defined.
\QED

\bigskip
\begin{defn}\label{defn}
As in Definition~\ref{defn3.4},
we define the invariant by
$$ GW^{J,{\cal H}}_{g,k}(X,A)(\,\beta;\alpha\,)\ =\
      [\,{\cal M}^{J,{\cal H}}_{g,k}(X,A)\,] \cap
      (\,\beta\cup \alpha\,).
$$
where $\alpha$ in $[\,H^{*}(X;{\Bbb Q})\,]^{k}$
and  $\beta$ in $H^{*}(\overline{{\cal M}}_{g,k};{\Bbb Q})$.
\end{defn}

\bigskip

These invariants are unchanged under  deformations
of complex structures which preserve the $(1,1)$-type of $A$.
In particular, if $A=mK$ for $m\geq 1$, then $GW^{J,{\cal H}}_{g,k}(X,A)$ is invariant under all deformations of the complex structure
since $mK$ is always of type $(1,1)$.

On the other hand, if $A$ is not a linear combination of components of
canonical class, then the set $F$ which appears in (\ref{frontier})
is empty, and hence
\begin{equation*}
 \CM_{g,k}^{J,\H}(X,A)\ =\ \Phi^{-1}(0)\ =\ \CM_{g,k}(X,A).
\end{equation*}
Therefore the  invariants
defined in Definitions~\ref{defn3.4} and ~\ref{defn}
coincide.  That completes the proof of Proposition~\ref{3}
stated in the introduction.

\bigskip
These invariants differ from the more familiar GW invariants.

\begin{ex}
The generic element of the linear system $|K|$ of the canonical class
is an embedded curve of genus $g=K^2+1$.   For that genus,  C. Taubes
proved that the (\,standard\,) GW invariant  $GW_{g}(X,K)$ is the same as
the Seiberg-Witten invariant and is equal to 1 for all surfaces of
general type.

On the other hand,  the family invariant for that genus vanishes.
That can be seen from a dimension count:  the linear system $|K|$ has
(real) dimension $2(p_g-1)$, while the family GW invariant
lies in dimension
$2\,(-K\cdot K +(g-1)+p_g)=2p_g$ (cf. Proposition~\ref{P:last}).
\end{ex}

\vskip 1cm


\setcounter{equation}{0} \setcounter{section}{4} \ {
\renewcommand{\theequation}{A.\arabic{equation}}
\renewcommand{\thetheorem}{A.\arabic{theorem}}
\section{Appendix -- Relations with the Behrend-Fantechi Approach}
\bigskip

   Behrend and Fantechi \cite{bf} have defined modified GW invariants
for K\"{a}hler surfaces using algebraic
geometry.  While their techniques are completely different from ours,
the definitions seem to be, at their
core, equivalent.  In this appendix we make several observations
which relate their approach to ours.  This is
necessarily tentative because the paper \cite{bf} is not yet available;
we are relying on the terse description
given in \cite{bl3}.

\medskip

In algebraic geometry, the virtual fundamental class
$[\overline{\M}_{g,k}(X,A)]^{\mbox{vir}}$
is obtained from the relative tangent-obstruction spaces
together with
the tangent-obstruction spaces of Deligne-Mumford space
$\overline{\M}_{g,k}$.
Behrend and Fantechi modified their machinery,
intrinsic normal cone and obstruction complex, by
replacing the relative obstruction space $H^{1}(f^{*}TX)$
by the kernel of  the map
\begin{equation}\label{BF-map}
H^{1}(f^{*}TX)\to H^{2}(X,{\cal O})
\end{equation}
   defined by  dualizing of the  composition
\begin{equation}\label{Dual}
H^{0}(X,\Omega^{2}) \to H^{0}(f^{*}\Omega^{2}) \to
H^{0}(f^{*}\Omega^{1}\otimes f^{*}\Omega^{1}) \to
H^{0}(f^{*}\Omega^{1}\otimes \Omega^{1}).
\end{equation}
In order for their machinery to work,
the map (\ref{BF-map}) is of constant rank --- in particular
surjective --- for every $f$ in
$\overline{\M}_{g,k}(X,A)$ \cite{bl3}.
Composing  (\ref{Dual}) with the Kodaira-Serre dual map, we have
\begin{equation}\label{K-map}
H^{0}(X,\Omega) \to
H^{0}(f^{*}\Omega^{1}\otimes \Omega^{1}) \to
H^{1}(f^{*}TX).
\end{equation}
This map is given by
$\beta\to K_{\beta}\,df\,j$.

\bigskip

\begin{proposition}
\label{A.1}
Let $(X,J)$ be a K\"{a}hler surface and $A\in H^{1,1}(X,\Z)$.
Then  the family moduli space
             $\overline{\M}^{\PH}_{g,k}(X,A)$ is compact if and only if
   the map (\ref{BF-map}) is surjective for every $f$ in
             $\overline{\M}_{g,k}(X,A)$.
\end{proposition}
\pf By Theorem~\ref{L:GT}   $\overline{\M}^{\PH}_{g,k}(X,A)$ consists
of pairs $(f,\alpha)$ with
$f\in \overline{\M}_{g,k}(X,A)$
and with the image of $f$ contained in  the zero set of $\alpha$; the
latter condition means that
$K_\alpha=0$  along the image, so $K_{\alpha}\,df\,j=0$ for all
$(f,\alpha)$.  As usual,
$\overline{\M}_{g,k}(X,A)$ is compact by the Gromov Compactness Theorem.

Now, suppose
(\ref{BF-map}) is surjective. Then by duality (\ref{K-map}) is
injective. This implies
$\alpha=0$ and hence
$\overline{\M}^{\PH}_{g,k}(X,A)=\overline{\M}_{g,k}(X,A)$ is compact.
Conversely, suppose for some
$f\in\overline{\M}_{g,k}(X,A)$ there is a $\beta$ in the kernel of
(\ref{K-map}).   Then setting
$\alpha=\beta+\bar\beta$ we have $\dbar f=tK_\alpha df j=0$ --- and
hence $(f,t\alpha)\in
\overline{\M}^{\PH}_{g,k}(X,A)$ --- for all real $t$.   That means
that  $\overline{\M}^{\PH}_{g,k}(X,A)$ is
compact only when  (\ref{K-map}) is injective or equivalently when
   (\ref{BF-map}) is surjective.
\QED

\bigskip

    The map (\ref{K-map}) is directly related to the linearization operator of
the $(J,\alpha)$-holomorphic map equation.

\medskip

   Suppose that $A$ is $(1,1)$ and that  the family moduli space
$\overline{\M}^{\PH}_{g,k}(X,A)$ is compact as in
Proposition \ref{A.1}. Consider the linearization of  the
$(J,\alpha)$-holomorphic map equation  at $(f,j,\alpha)$.
Since $J$ is K\"{a}hler, the linearization reduces to
$$
L_{f}\oplus Jdf \oplus L_{0} :
         \Omega^{0}(f^{*}TX)\oplus
          T_{j}\overline{\cal M}_{g,n}\oplus \PH \to
                      \Omega^{0,1}(f^{*}TX)
\mbox{\ \ where}\ \
\left\{\begin{array}{lll}
   L_{f}(\xi) & = & \nabla\xi + J\nabla\xi\, j \\ \\
    L_{0}(\beta)  & = &  -2K_{\beta}df j
\end{array}\right.
$$
In fact, this $L_{f}$ is exactly (twice) the Dolbeault derivative
$\overline{\partial}$. Therefore,
$\mbox{Ker}(L_{f})$ are $\mbox{Coker}(L_{f})$ are identified with  the
Dolbeault cohomology groups $H^{0} (f^{*}TX)$
and $H^{0,1}(f^{*}TX)$, respectively.

\begin{proposition}
\label{A.2}
Under either of the  two equivalent conditions of Proposition
\ref{A.1}  there are natural identifications
$H^{1}(f^{*}TX) \simeq H^{0,1}(f^{*}TX)
$ and $H^{0}(X,\Omega^{2}) \simeq \PH$
    under which identification the map is identified with $(\ref{K-map})$ with
$$
L_{0}:\PH\to \mbox{Coker}(L_{f}\oplus Jdf).
$$
By Proposition \ref{A.1} this map is injective if and only if the
family moduli space
             $\overline{\M}^{\PH}_{g,k}(X,A)$ is compact.
\end{proposition}

\pf
It follows by comparing the formulas for $L_0$ and
$(\ref{K-map})$ that $L_{0}$ maps  $\PH$ into
$\mbox{Coker}(L_{f})$. On the other hand,
given $h\in T_{j}\overline{\cal M}_{g,n}$,
there is a family $j_{t}$ with $j_{0}=j$ and
$\frac{d\,j_{t}}{d\,t}_{|_{t=0}}=h$.
It follows from Proposition~\ref{P:Basic}b
and $\langle \beta, A\rangle = 0$ that
\begin{equation*}
0= \left. \frac{d}{d\,t}\right |_{t=0}
     \int_{(\Sigma,j_{t})}f^{*}(\beta)
=  \left.\frac{d}{d\,t}\right |_{t=0}
     \int_{(\Sigma,j_{t})}
     \langle df + Jdfj_{t}  , K_{\beta}{f}_{*} j_{t} \rangle
=  \int_{\Sigma}
     \langle Jdf(h),K_{\beta}f_{*} j \rangle.
\end{equation*}
This implies that $L_{0}$ maps  $\PH$ into
$\mbox{Coker}(L_{f}\oplus Jdf)$.
\QED



\begin{thebibliography}{}







\bibitem[B]{b} L. Besse, {\em Einstein manifolds},
                  Springer-Verlag, Berlin Heidelberg, 1987.

\bibitem[BL1]{bl1} J. Bryan and N. C. Leung,
                      {\em The enumerative geometry of K3 surfaces
                           and modular forms}, J. Amer. Math. Soc.
                      {\bf 13} (2000), 371-410.

\bibitem[BL2]{bl2} J. Bryan and N. C. Leung, {\em Genenerating functions
for the number of curves on abelian surfaces},
Duke Math. J. {\bf 99} (1999), no. 2, 311-328.

\bibitem[BL3]{bl3} J. Bryan and N. C. Leung, {\em Counting curves on
irrational surfaces}, Surveys in differential geometry:
differential geometry inspired by string theory, 313-339,
Surv. Diff. Geom., 5, Int. Press, Boston, MA, 1999.



\bibitem[BF]{bf} K. Behrend and B. Fantechi, In Preparation.

\bibitem[BPV]{bpv} W. Barth, C. Peters, and A. Van de Ven,
                      {\em Compact complex surfaces},
                       Springer-Verlag, Berlin Heidelberg, 1984.



\bibitem[D]{d} S. K. Donaldson,
{\em Yang-Mills invariants of 4-manifolds },
In S.K. Donaldson and C.B. Thomas, editors,
Geometry of low dimensional manifolds : Gauge Theory and Algebraic surfaces,
number 150 in London Math. Soc. Lecture Note Series,
Cambrige University Press, 1989.




\bibitem[FM]{fm} R. Friedman, J.W. Morgan, {\em Smooth four-manifolds
       and complex surfaces},
       Springer-Verlag, Berlin Heidelberg, 1994.






\bibitem[G]{g} L. G\"{o}ttsche, {\em A conjectural generating
function for numbers of curves on surfaces}, preprint,
alg-geom/9711012




\bibitem[GH]{gh} P. Griffiths and J. Harris, {\em Principles of algebraic
geometry}, J. Wiley, New York, 1978.






\bibitem[IP1]{ip1} E. Ionel and T.  Parker,  {\em Relative
Gromov-Witten Invariants}, Ann. Math. {\bf 157} (2003), 45-96


\bibitem[IP2]{ip2} E. Ionel and T.  Parker,  {\em The Symplectic Sum
Formula for Gromov-Witten Invariants}, preprint, math.SG/0010217.


\bibitem[IS]{is} S. Ivashkovich and V. Shevchishin,
{\em Gromov compactness theorem for $J$-complex curves with boundary},
Internat. Math. Res. Notices {\bf 2000}, no. 22, 1167-1206.


\bibitem[KM]{km} P. Kronheimer and T. Mrowka, {\em Embedded surfaces and the
structure of
Donaldson's polynomial invariants}, J. Differential Geom.
{\bf 41} (1995), 573-734.



\bibitem[KP]{kp} S. Kleiman and R. Piene,
{\em Enumerating singular curves on surfaces},
Algebraic Geometry : Hirzebruch 70 (Warsaw, 1998), 209-238,
Contemp. Math., 241, Amer. Math. Soc., Providence, RI, 1999.









\bibitem[LT]{LT} J. Li and G. Tian,  {\em Virtual moduli cycles and
Gromov-Witten invariants of general symplectic manifolds}, Topics in
symplectic $4$-manifolds (Irvine, CA, 1996), 47--83, First Int.
Press Lect. Ser., I, International Press, Cambridge, MA, 1998.







\bibitem[P]{p} T. Parker, {\em Compactified moduli spaces of
pseudo-holomorphic curves} Mirror symmetry, III (Montreal,
PQ, 1995), 77--113, AMS/IP Stud. Adv. Math., 10, Amer. Math. Soc.,
Providence, RI, 1999.

\bibitem[PW]{pw} T. Parker and J. Wolfson,
{\em Pseudo-holomorphic maps and bubble trees},
Jour. Geometric Analysis, {\bf 3} (1993) 63-98.









\bibitem[V]{v} I. Vainsencher, {\em Enumeration of n-fold tangent
hyperplanes to a surface}, J. Alg. Geom. {\bf 4} (1995), 503-526


\bibitem[YZ]{yz} S.T. Yau and E. Zaslow, {\em BPS States, String Duality,
and Nodal Curves on K3}, Nuclear Phys. B {\bf 471} (1996), 503-512




\end{thebibliography}
\end{document}